\documentclass[11pt]{article} 
\usepackage[utf8]{inputenc}
\usepackage{caption,graphicx,amsmath,amsthm,amssymb, authblk,xcolor,enumitem,subfigure,float,verbatim,hyperref,pgfplots,bibentry,mathtools,soul,isomath, cancel,multicol,wrapfig,subcaption}
\usepackage{mathrsfs}
\usepackage{indentfirst}
\usepackage{nccmath}
\graphicspath{ {images/} }
\usepackage{multicol}
\usepackage{graphicx,amsmath,amssymb,amsthm,marvosym,tikz,fontenc}
\usetikzlibrary{backgrounds}
\usetikzlibrary{patterns,shapes.misc, positioning}
\usepackage{pgfplots}
\pgfplotsset{width=10cm,compat=1.9,tick scale binop=\times}
\usepackage{indentfirst}
\usepackage[a4paper,bindingoffset=0.2in,%
left=0.8in,right=1in,top=1.2in,bottom=1.2in,%
footskip=.25in]{geometry}
\usepackage{relsize}
\usepackage[english]{babel}
\allowdisplaybreaks
\theoremstyle{plain}
\newtheorem{theorem}{Theorem}[section]
\newtheorem{lemma}[theorem]{Lemma}

\newtheorem{corollary}[theorem]{Corollary}
\newtheorem{definition}[theorem]{Definition}

\newtheorem{example}[theorem]{Example}
\DeclareRobustCommand{\rchi}{{\mathpalette\irchi\relax}} 
\newcommand{\irchi}[2]{\raisebox{\depth}{$#1\chi$}}   

\newlength{\defbaselineskip}
\newcommand{\setlinespacing}[1]%
{\setlength{\baselineskip}{#1 \defbaselineskip}}

\date{}
\begin{document}
\title{On the Construction of Singular and Cospectral Hypergraphs}
\author{Liya Jess Kurian$^1$\footnote{liyajess@gmail.com},  Chithra A. V$^1$\footnote{chithra@nitc.ac.in}
 \\ \small 
 1 Department of Mathematics, National Institute of Technology Calicut,\\\small
 Calicut-673 601, Kerala, India\\ \small	}
\maketitle
\thispagestyle{empty}
 \begin{abstract}
 In this paper, we define two operations, neighbourhood m-splitting hypergraph $NS_m(\mathscr{G}^*)$ and non-neighbourhood splitting hypergraph $NNS(\mathscr{G}^*)$, and obtain several properties of their adjacency spectrum. We also estimate the energies of $NS_m(\mathscr{G}^*)$ and $NNS(\mathscr{G}^*)$. Moreover, we introduce two new join operations on $k$-uniform hypergraphs: the neighbourhood splitting V-vertex join $\mathscr{G}_1^*\veebar \mathscr{G}_2^*$ and the S-vertex join $\mathscr{G}_1^*\barwedge \mathscr{G}_2^*$ of hypergraphs $\mathscr{G}_1^*$ and $\mathscr{G}_2^*$, and determine their adjacency spectrum. As an application, we obtain infinite families of singular hypergraphs and infinite pairs of non-regular non-isomorphic cospectral hypergraphs.\\
  
\end{abstract}
\section{Introduction}
Let $\mathscr{G}^*=(V,E)$ be a hypergraph of order $n$, where $V(\mathscr{G}^*)=\{v_1,v_2,v_3,\ldots,v_n\}$ is the vertex set  and the collection of hyperedges $E(\mathscr{G}^*)=\{e_1,e_2,e_3,\ldots,e_t\}$  is the edge set of $\mathscr{G}^*$. Each hyperedge in $E(\mathscr{G}^*)(|e_i| \geq 2)$ is a non-empty subset of the vertex set $V(\mathscr{G}^*)$. Throughout this paper, we consider $k$-uniform hypergraphs, which are hypergraphs whose hyperedges contain exactly $k(k\geq2)$ vertices\textnormal{\cite{Cooper2012,Kumar2017}}. When $k=2$ it becomes an ordinary graph. The degree of a vertex $v\in V$, $d(v)$, is defined as the number of hyperedges which contain the vertex $v$. A hypergraph in which every vertex $v_i\in V$  has degree $r$ is said to be a $r$-regular hypergraph. If a hypergraph is both $k$-uniform and $r$-regular, we refer to it as a $(k,r)$-regular hypergraph. In \cite{Kumar2017}, the authors focus on the characteristics of $(k,r)$ regular hypergraphs. A $k$-uniform hypergraph $\mathscr{G}^*$ with $n$ vertices is said to be a complete $k$-uniform hypergraph $K_n^k$ if $E(\mathscr{G}^*)$ is the collection of all possible $k$-subsets of $V(\mathscr{G}^*)$ \cite{Berge1973}.  Let $D$ be any $(k-1)$-subset of $V(\mathscr{G}^*)$ of hypergraph $\mathscr{G}^*$. Then the vertex $v\in V$ is said to be the neighbour of $D$ ($vND$) if $\{v,D\} \in E,$ and otherwise it is not a neighbour of $D$ ($v\cancel N D$).\\

 The adjacency matrix \cite{Bretto2013}  of $\mathscr{G}^*$, $A(\mathscr{G}^*)$, is a square matrix of order $n$ whose rows and columns are indexed by the vertices of $\mathscr{G}^*$. For all $ v_i,v_j \in V, $
 \begin{equation*}
a_{ij} =\left\{
   \begin{array}{ll}
      \mid \{e_k  \in E(\mathscr{G}^*): \{v_i,v_j\} \subset  e_k\}\mid    &  \mbox{, } v_i \neq v_j, k\in[1,t]  \\
       0  & \mbox{, } v_i=v_j
   \end{array}.
   \right.
\end{equation*}
Clearly, it generalises the definition of the adjacency matrix of graphs. A scalar $\lambda$ is an eigenvalue of a matrix $M$ if there exists a non-zero eigenvector $\mathbf{x}$ that satisfies the equation $M\mathbf{x}=\lambda \mathbf{x} $. The adjacency spectrum refers to the collection of all eigenvalues of the matrix $A(\mathscr{G}^*)$ along with their corresponding multiplicities. Let $\lambda_{1},\lambda_{2},\lambda_{3},\ldots,\lambda_{d}(d\leq n)$ be the eigenvalues, and $m_{1},m_{2},m_{3},...,m_{d}$ be the corresponding multiplicities of the adjacency matrix $A(\mathscr{G}^*)$. Then adjacency spectrum of $\mathscr{G}^*$, $\sigma_{A}(\mathscr{G}^*)$ is given by,
$$ \sigma_{A}(\mathscr{G}^*)=\begin{pmatrix}
               \lambda_{1} & \lambda_{2} & \lambda_{3} & \cdots & \lambda_{d}\\
               m_{1} & m_{2} & m_{3} & \cdots  & m_{d}
             \end{pmatrix}.$$
If the spectrum of hypergraphs $\mathscr{G}^*$ and $\mathcal{H}^*$ coincide, they are cospectral. The spectral radius of $\mathscr{G}^*\,(\rho(A(\mathscr{G}^*)))$ is the largest absolute value of the eigenvalues of $A(\mathscr{G}^*)$. The energy(adjacency energy) $\mathcal{E}(\mathscr{G}^*)$ of $\mathscr{G}^*$ is defined as the sum of the absolute values of the adjacency eigenvalues of $\mathscr{G}^*$. A hypergraph $\mathscr{G}^*$ is singular if it has zero as an adjacency eigenvalue and the multiplicity of 0 is the nullity of $\mathscr{G}^*,\,\eta(\mathscr{G}^*)$. In \cite{Gutman2016}, the authors have conjectured that for a graph, the energy decreases when nullity increases. The study of singular graphs is a significant mathematical problem in relation to molecular orbital theory and network theory\cite{Tarimshawy2018,Sciriha2007,Sciriha2014}.\\

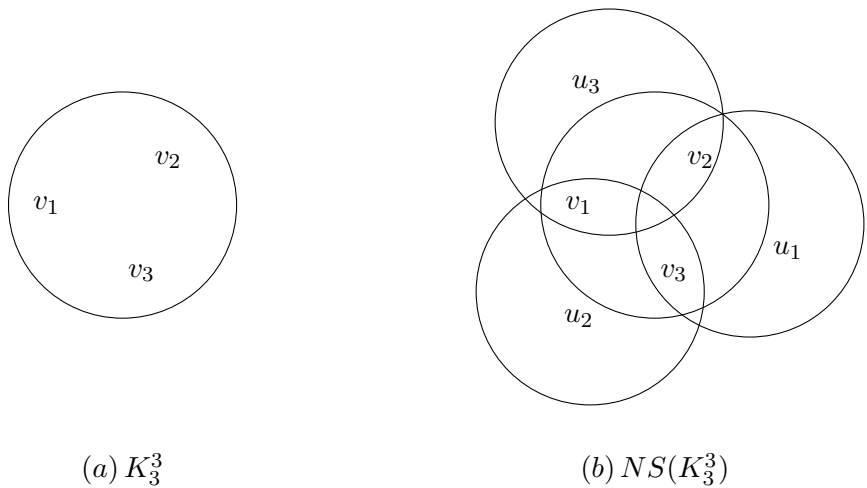
\begin{figure}[H]
\centering
\begin{tikzpicture}[scale=0.5]
     \draw (0,0) circle (3);
     \node at (1.2,1.2) {$v_2$};
     \node at (-2,0) {$v_1$};
     \node at (0.5,-1.8) {$v_3$};
    \draw (-1.2,2.2) circle (3);
    \node at (-1.8,3.2) {$u_3$};
    \draw (2.5,-0.5) circle (3);
    \node at (3.5,-1.2) {$u_1$};
     \draw (-1.7,-2.3) circle (3);
     \node at (-2,-3) {$u_2$};
 \node at (-0,-7) {$(b)\,NS(K^3_3)$};
\draw (-14,0) circle (3);
     \node at (-12.8,1.2) {$v_2$};
     \node at (-16,0) {$v_1$};
     \node at (-13.5,-1.8) {$v_3$};
 \node at (-14,-7) {$(a)\,K_3^3$};
\end{tikzpicture}
  \caption{Neighbourhood splitting hypergraph of $K_3^3$  }
  \label{fig:1}
\end{figure}

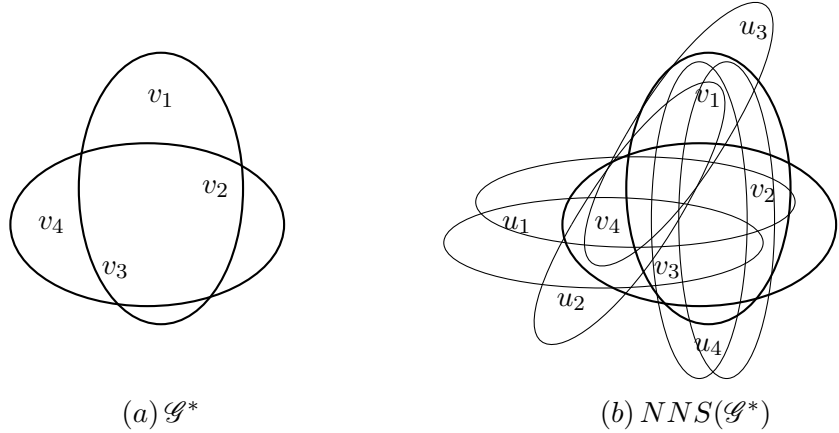
\begin{figure}[H]
\centering
\begin{tikzpicture}[scale=0.6]
   \draw[thick] (0,0) ellipse (1.8 and 3);
     \node at (1.2,0) {$v_2$};
     \node at (-1.0,-1.8) {$v_3$};
     \node at (0,2) {$v_1$};
     \draw[thick] (-0.3,-0.8) ellipse (3 and 1.8);
      \node at (-2.4,-0.8) {$v_4$};
\node at (0,-5) {$(a)\,\mathscr{G}^*$};
      \draw[thick] (12,0) ellipse (1.8 and 3);
     \node at (13.2,-0.1) {$v_2$};
     \node at (11.1,-1.8) {$v_3$};
     \node at (12,2) {$v_1$};
     \draw[thick] (11.8,-0.8) ellipse (3 and 1.8);
      \node at (9.8,-0.8) {$v_4$};

       \draw (9.7,-1.2) ellipse (3.5 and 1);
       \draw (10.4,-0.3) ellipse (3.5 and 1);
       \node at (7.8,-0.8) {$u_1$};
       \draw(11.8,-0.7) ellipse (1.05 and 3.5);
       \draw (12.4,-0.7) ellipse (1.05 and 3.5);
       \node at (12,-3.5) {$u_4$};
       \node at (13,3.5) {$u_3$};
       \draw (10.2,1.2) to [out=60,in=60, looseness=5] (12.5,1.2);
       \draw (10.2,1.2) to [out=240,in=240, looseness=5](12.5,1.2);
       \draw (9.,-0.8) to [out=60,in=60, looseness=5] (11.4,-0.6);
       \draw (9.,-0.8) to [out=240,in=240, looseness=4.5](11.4,-0.6);
       \node at (9,-2.5) {$u_2$};
       \node at (11.5,-5) {$(b)\,NNS(\mathscr{G}^*)$};
\end{tikzpicture}

  \caption{Non-neighbourhood splitting hypergraph of $\mathscr{G}^*$   }
  \label{fig:2}
\end{figure}

Recently, researchers have shown great interest in determining the cospectral family of hypergraphs. Neighbourhood splitting V-vertex and S-vertex join, non-neighbourhood splitting vertex join, central vertex join, and central edge join are some of the operations defined in graph theory to find the cospectral family of graphs\cite{Berman2024,jahfar2020,Lu2023}. Motivated by these researches, we found new cospectral families of hypergraphs by extending the study of neighbourhood splitting V-vertex and S-vertex join of graphs. To study the spectral properties of hypergraphs, we introduced two new graph operations neighbourhood m-splitting hypergraphs and non-neighbourhood splitting hypergraphs. 
\begin{definition}
    Let $\mathscr{G}^*$ be a $k-$uniform hypergraph with vertex set $V(\mathscr{G}^*)=\{v_1,v_2,v_3,\ldots,v_n\}$ and edge set $E(\mathscr{G}^*)$. The neighbourhood splitting hypergraph $NS(\mathscr{G}^*)$ is a hypergraph with vertex set $V(\mathscr{G}^*)\cup V'$ and edge set $E(\mathscr{G}^*)\cup E'$, where $V'=\{u_1,u_2,u_3,\ldots,u_n\}$ and $E'=\{\{u_i,D\}: v_i N D, u_i\in V', D\subset V  \}$.
\end{definition}
The neighbourhood splitting hypergraph is obtained by introducing a new vertex corresponding to each vertex of $\mathscr{G}^*$. The collection of all such new vertices is denoted by $S(\mathscr{G}^*)$. For example, the neighbourhood splitting hypergraph of $K^3_3$ is given in Figure \ref{fig:1}.
\begin{definition}
    Let $\mathscr{G}^*$ be a $k-$uniform hypergraph with vertex set $V=\{v_1,v_2,v_3,\ldots,v_n\}$ and edge set $E$. The non-neighbourhood splitting hypergraph $NNS(\mathscr{G}^*)$ is a hypergraph with vertex set $V\cup V'$ and edge set $E\cup E'$, where $V'=\{u_1,u_2,u_3,\ldots,u_n\}$ and $E'=\{\{u_i,D\}: v_i \cancel N D, u_i\in V', D\subset V  \}$.
\end{definition}
\noindent The non-neighbourhood splitting hypergraph of $\mathscr{G}^*$ is shown in Figure \ref{fig:2}.
This paper aims to analyse the energy of hypergraphs obtained by extending the spectral theory of splitting graphs to hypergraphs and to construct a family of singular hypergraphs and cospectral hypergraphs using new operations defined on neighbourhood splitting hypergraphs. We organise the remaining sections of this paper as follows: In Section 2, we provide all the necessary definitions and lemmas. In Section 3, we focus on analysing the spectrum of neighbourhood m-splitting and non-neighbourhood splitting hypergraphs and evaluating their energy.  In addition, we give results to obtain infinite families of singular hypergraphs based on the spectrum of neighbourhood m-splitting hypergraphs. In Section 4, we introduce the operations neighbourhood splitting (V-vertex and S-vertex) join and determine their eigenvalues. Also, we construct pairs of non-isomorphic cospectral hypergraphs. Further, we give a result on the singularity of the neighbourhood splitting join of hypergraphs. 
\section{Preliminaries}
\noindent \textbf{Notation:} We write $j\in [a,b]$ if $j$ takes all the integer values satisfying the condition $a\leq j\leq b$.
Let $J_n$ and $I_n$ denote the all one and identity matrix of order $n$ and $ J_{k,n}$ denote all one matrix of order $k\times n$, respectively. 
 \begin{lemma}\label{detblock}\textnormal{\cite{Das2018}}
Let $M_{11}, M_{12}, M_{21}$, and $M_{22}$ be matrices with $M_{11}$ invertible. Let 
             $$ M= \begin{bmatrix}
                 M_{11} & M_{12}\\
                M_{21} & M_{22}
                  \end{bmatrix}$$
Then, $\det(M)=det(M_{11})det(M_{22}-M_{21}M_{11}^{-1}M_{12})$ and if $M_{22}$ is invertible, then $det(M) = det(M_{22})det(M_{11}-M_{12}M_{22}^{-1}M_{21}).$
          \end{lemma}
\begin{definition}\textnormal{\cite{McLeman2011}}\label{coronaldefn}
    The coronal $\rchi_{M}(\lambda)$ of a $n\times n$ matrix $M$ is defined as the sum of the matrix $(\lambda I_n -M)^{-1}$, that is,
         $$\rchi_{M}(\lambda)=J_{1,n}(\lambda I-M)^{-1}J_{n,1}.$$

\end{definition}
\begin{lemma}\textnormal{\cite{McLeman2011}}\label{coronalval}
    The coronal $\rchi_{M}(\lambda)$ of a $n\times n$ matrix $M$, whose row sum is equal to $r$. Then
    $$\rchi_{M}(\lambda)=\frac{n}{\lambda-r}.$$
\end{lemma}
\begin{lemma}\textnormal{\cite{Cvetkovic2010}}\label{detMJ}
    Let M be an $n\times n$ real matrix and $\alpha\in \mathbb{R}$. Then
    $$\det(M+\alpha J_n)=\det(M)+ \alpha J_{1,n} adj(M) J_{n,1}.$$
\end{lemma}
\begin{lemma}\textnormal{\cite{jahfar2020}}\label{invIJ}
              For any two real numbers $r$ and $s$,
              $$(rI_n-sJ_n)^{-1}=\frac{1}{r}I_n+\frac{s}{r(r-ns)}J_n.$$
\end{lemma}
\begin{definition}\textnormal{\cite{Horn1994}}
     Let $P=(p_{ij})$ and $Q$ be two matrices of any order. Then the Kronecker product of $P$ and $ Q$ is a block matrix,
$$P\otimes Q=(p_{ij} B).$$
If $\lambda$ and  $\beta$ are the eigenvalues of $P$ and $Q$ respectively, then $\lambda\beta$ is an eigenvalue of $P\otimes Q$.   
\end{definition}

\section{Neighbourhood m-Splitting and Non-neighbourhood Splitting Hypergraphs}
This section determines the spectrum of neighbourhood and non-neighbourhood splitting hypergraphs and computes their adjacency energy. 
\begin{definition}
    Let $\mathscr{G}^*$ be a $k-$uniform hypergraph with vertex set $V=\{v_1,v_2,v_3,\ldots,v_n\}$ and edge set $E$. The neighbourhood m-splitting hypergraph $NS_m(\mathscr{G}^*)$ is a hypergraph with vertex set $V\cup\Bigl(\cup_{i=1}^m V'_i\Bigr)$ and edge set $E\cup\Bigl(\cup_{i=1}^m E'_i\Bigr)$, where $m\geq 1$, and $V'_i=\{u_{i1},u_{i2},u_{i3},\ldots,u_{in}\}$ and $E'_i=\{\{u_{ij},D\}: v_j N D, u_{ij}\in V'_i, D\subset V ,j\in[1,n] \}$.
\end{definition}
The adjacency matrix of $NS_m(\mathscr{G}^*)$ can be represented as a block matrix, based on the ordering of the vertices $V(\mathscr{G}^*)$ and $V'_i,i\in[1,n]$, as follows:
\begin{equation*}
      A(NS_m(\mathscr{G}^*))=\begin{bmatrix}
	(mk-2m+1)A(\mathscr{G}^*) & A(\mathscr{G}^*)& A(\mathscr{G}^*)&\cdots& A(\mathscr{G}^*)\\
        A(\mathscr{G}^*) & \mathbf{0}_n& \mathbf{0}_n&\cdots& \mathbf{0}_n\\
         A(\mathscr{G}^*) & \mathbf{0}_n& \mathbf{0}_n&\cdots& \mathbf{0}_n\\
        \vdots & \vdots& \vdots&\ddots& \vdots\\
	A(\mathscr{G}^*) & \mathbf{0}_n&\mathbf{0}_n&\cdots&\mathbf{0}_n
\end{bmatrix},
\end{equation*}
where $\mathbf{0}_n$ is the zero matrix of order $n$. 
\begin{theorem}\label{SpecNmS}
    Let $\mathscr{G}^*$ be a $k$-uniform hypergraph of order $n$ with eigenvalues $\lambda_i,\,i\in[1,n]$. Then the spectrum of $NS_m(\mathscr{G}^*)$ is given by,
    $$ \sigma_{A}(NS_m(\mathscr{G}^*))=\begin{pmatrix}
              0&\frac{mk-2m+1+\sqrt{(mk-2m+1)^2+4m}}{2}\lambda_i& \frac{mk-2m+1-\sqrt{(mk-2m+1)^2+4m}}{2}\lambda_i\\
               n(m-1) & 1 & 1
             \end{pmatrix}.$$
\end{theorem}
 \begin{proof}
   Let $M=\begin{bmatrix}
         mk-2m+1 & 1 & 1 & 1& \cdots&1\\
         1 & 0 & 0 & 0& \cdots&0\\
          1 & 0 & 0 & 0& \cdots&0\\
           \vdots & \vdots & \vdots &\vdots& \ddots&\vdots\\
            1 & 0 & 0 & 0& \cdots&0
     \end{bmatrix}$.  Then the adjacency matrix of $NS_m(\mathscr{G}^*)$ can be written as,
     $$A(NS_m(\mathscr{G}^*)=M\otimes A(\mathscr{G}^*).$$
      Since rank$(M)=2$ and trace$(M)=mk-2m+1$, then sum of the non-zero eigenvalues $\nu_1$ and $\nu_2$ of $M$ is 
     \begin{equation}\label{eqn7.1}
         \nu_1+\nu_2=mk-2m+1.
     \end{equation}
  Again,
  $$M^2=\begin{bmatrix}
         (mk-2m+1)^2+m  & mk-2m+1 & mk-2m+1& \cdots&mk-2m+1\\
         mk-2m+1& 1 & 1& \cdots&1\\
          mk-2m+1  & 1 & 1& \cdots&1\\
           \vdots  & \vdots &\vdots& \ddots&\vdots\\
            mk-2m+1  & 1 & 1& \cdots&1
     \end{bmatrix},$$
and we have 
\begin{equation}\label{eqn8}
   trace(M^2)=\nu_1^2+\nu_2^2=(mk-2m+1)^2+2m. 
\end{equation}
From (\ref{eqn7.1}) and (\ref{eqn8}), we get
$$\nu_1=\frac{mk-2m+1+\sqrt{(mk-2m+1)^2+4m}}{2},\:\;\nu_2=\frac{mk-2m+1-\sqrt{(mk-2m+1)^2+4m}}{2},$$
and all the other $m-1$ eigenvalues of $M$ are zeros. If $\lambda_i,\,i\in[1,n]$ are the eigenvalues of $A(\mathscr{G}^*)$, then $\lambda_i\nu_1, \lambda_i\nu_2$ and 0 with multiplicity $n(m-1)$ are the eigenvalues of $M\otimes A(\mathscr{G}^*).$ Thus the theorem follows.
\end{proof}
 In \cite{Vaidya2017}, authors established the spectrum and energy of $m$-splitting graphs. We denote the neighbourhood 1-splitting hypergraph (neighbourhood splitting hypergraph) by $NS(\mathscr{G}^*)$. The proof of the Corollary \ref{specNS} is obtained from Thereom \ref{SpecNmS} for $m=1.$
\begin{corollary}\label{specNS}
    Let $\mathscr{G}^*$ is a $k$-uniform hypergraph of order $n$ with eigenvalues $\lambda_i,\,i\in[1,n]$. Then the adjacency spectrum of $NS(\mathscr{G}^*)$ is given by
    \begin{equation*}
        \sigma_{A}(NS(\mathscr{G}^*))=\begin{pmatrix}
              \hat{h}\lambda_i  & \hat{\hat{h}}\lambda_i\\
               1 & 1
             \end{pmatrix},
    \end{equation*} 
             where $\hat{h}=\frac{(k-1)+\sqrt{(k-1)^2+4}}{2}$ and $\hat{\hat{h}}=\frac{(k-1)-\sqrt{(k-1)^2+4}}{2}.$
\end{corollary}

Next, we determine the energy of $NS_m(\mathscr{G}^*)$ in terms of the energy of $\mathscr{G}^*$.
\begin{corollary}\label{EnergyNS}
    Let $\mathscr{G}^*$ be a $k$-uniform hypergraph of order $n$ with eigenvalues $\lambda_i,\,i\in[1,n]$. Then the energy of $NS_m(\mathscr{G}^*)$ is given by,
    $$\mathcal{E}(NS_m(\mathscr{G}^*))=\sqrt{(mk-2m+1)^2+4m}\,\mathcal{E}(\mathscr{G}^*).$$
\end{corollary}
\begin{proof}
    \begin{align*}
        \mathcal{E}(NS_m(\mathscr{G}^*))&=\Bigl(\Bigl|\frac{mk-2m+1+\sqrt{(mk-2m+1)^2+4m}}{2}\Bigr|\\
        &\hspace{4cm}+\Bigl|\frac{mk-2m+1-\sqrt{(mk-2m+1)^2+4m}}{2}\Bigr|\Bigr)\sum_{i=1}^n|\lambda_i|\\
        &=\sqrt{(mk-2m+1)^2+4m} \,\mathcal{E}(\mathscr{G}^*).
    \end{align*}
\end{proof}
From Theorem \ref{SpecNmS}, it is clear that the nullity of $NS_m(G^*)$  increases with $m$. 
\begin{corollary}
    Let $\mathscr{G}^*$ be a $k$-uniform hypergraph. Then $\eta(\mathscr{G}^*)\leq \eta(NS_{m}(\mathscr{G}^*))$, when $m=1$ equalty holds.
\end{corollary}
Now, we investigate the properties of the eigenvalues of $NS_m(\mathscr{G}^*)$. Note that  $NS_m(\mathscr{G}^*)$ has $n(m-1)$ eigenvalues that are equal to zero, and the remaining $2n$ eigenvalues of $NS_m(\mathscr{G}^*)$ possess certain properties.
\begin{theorem}\label{thm3.3}
  Let $\mathscr{G}^*$ be a $k$-uniform hypergraph of order $n$ with eigenvalues $\lambda_i,\,i\in[1,n]$. Then,
  \begin{enumerate}
      \item If $\lambda$ is any non-zero adjacency eigenvalue of $NS_m(\mathscr{G}^*)$ then there exist $\lambda_i$ such that $\displaystyle\frac{-m\lambda_i^2}{\lambda}$ is an adjacency eigenvalue of $NS_m(\mathscr{G}^*)$.
      \item $\lambda$ is an adjacency eigenvalue of $NS_m(\mathscr{G}^*)$ if and only if $-\displaystyle\frac{\nu_2^2}{m}\lambda$ is an eigenvalue of $NS_m(\mathscr{G}^*)$.
  \end{enumerate}
\end{theorem}
\begin{proof}
    From Theorem \ref{SpecNmS}, we get all the eigenvalues of $NS_m(\mathscr{G}^*)$. It can be noted that 
    \begin{align*}
        \frac{(mk-2m+1+\sqrt{(mk-2m+1)^2+4m})\lambda_i}{2}\,.&\frac{(mk-2m+1-\sqrt{(mk-2m+1)^2+4m})\lambda_i}{2}\\
        &\hspace{5.9cm}=-m\lambda_i^2.
    \end{align*}
    and hence the first statement holds. Let $\lambda$ be any non-zero eigenvalue of $NS_m(\mathscr{G}^*)$ (say, $\lambda=\frac{(mk-2m+1+\sqrt{(mk-2m+1)^2+4m})\lambda_i}{2}$). Since
    $$\frac{\frac{(mk-2m+1-\sqrt{(mk-2m+1)^2+4m})\lambda_i}{2}}{\frac{(mk-2m+1+\sqrt{(mk-2m+1)^2+4m})\lambda_i}{2}}=\displaystyle\frac{(mk-2m+1-\sqrt{(mk-2m+1)^2+4m})^2}{-4m}= \frac{-\nu_2^2}{m},$$
    then $-\displaystyle\frac{\nu_2^2}{m}\lambda$ is also an adjacency eigenvalue of $NS(\mathscr{G}^*)$.
\end{proof}
The problem of characterising singular graphs is challenging. For the last eighty years, scholars have discussed the structure of singular graphs\cite{Sciriha2007} and their construction\cite{Sciriha1998}. A family of singular hypergraphs can be obtained from neighbourhood m-splitting hypergraphs. The following corollary is the direct consequence of the  Theorems \ref{SpecNmS} and \ref{thm3.3}. Here, we discuss the relation between the spectral radius of $\mathscr{G}^*$ and $NS_m(\mathscr{G}^*)$ and the singularity of $NS_m(\mathscr{G}^*)$ as an extension of the study of singular graphs. 
\begin{corollary}
    Let $\mathscr{G}^*$ be a $k$-uniform hypergraph of order $n$ with eigenvalues $\lambda_i,\,i\in[1,n]$. Then
    \begin{enumerate}
        \item $\det\left(A(NS(\mathscr{G}^*))\right)=(-1)^n\det(A(\mathscr{G}^*))^2$ and singularity of $NS(\mathscr{G}^*)$ depends upon the singularity of $\mathscr{G}^*$
        \item  When $m>1$, $\det\left(A(NS_m(\mathscr{G}^*))\right)=0$ and is always singular .
        \item Adjacency-spectral radius of $NS_m(\mathscr{G}^*)$ depends on adjacency-spectral radius of $\mathscr{G}^*$, $$\rho(NS_m(\mathscr{G}^*))=\frac{mk-2m+1+\sqrt{(mk-2m+1)^2+4m}}{2} \rho(\mathscr{G}^*) .$$
    \end{enumerate}
\end{corollary}
Next, we discuss about the adjacency matrix and spectrum of $NNS(\mathscr{G}^*)$. The adjacency matrix of $NNS(\mathscr{G}^*)$, for $k\geq 3$ can be represented as a block matrix as follows:
\begin{equation*}
      A(NNS(\mathscr{G}^*))=\begin{bmatrix}
	(n-2)\binom{n-3}{k-3}(J_n-I_n)-(k-3)A(\mathscr{G}^*) & \binom{n-2}{k-2}(J_n-I_n)-A(\mathscr{G}^*)\\
	\binom{n-2}{k-2}(J_n-I_n)-A(\mathscr{G}^*) & \mathbf{0}_n
\end{bmatrix},
\end{equation*}
where $J_n$ and $\mathbf{0}_n$ denotes all one matrix and the zero matrix of order $n$, respectively.
\begin{theorem}
    Let $\mathscr{G}^*$ be a $(k,r)-$regular hypergraph of order $n$ with eigenvalues $\lambda_1=r(k-1)\geq \lambda_2\geq \lambda_3\geq \ldots \geq \lambda_n$. For $k\geq 3$ the adjacency spectrum of $NNS(\mathscr{G}^*)$ is given by
    $$ \sigma_{A}(NNS(\mathscr{G}^*))=\begin{pmatrix}
             \frac{-\bigl((k-3)\lambda_i+(n-2)\binom{n-3}{k-3}\bigr)\pm\sqrt{\bigl((k-3)\lambda_i+(n-2)\binom{n-3}{k-3}\bigr)^2+4\bigl(\lambda_i+\binom{n-2}{k-2}\bigr)}\bigr)^2}{2}  & \alpha_1&\alpha_2\\
               1& 1& 1
             \end{pmatrix},$$
             where $i\in [2,n]$, and $\alpha_1$ and $\alpha_2$ are the roots of the equation $\lambda^2+\Bigl(r(k-1)(k-3)-(n-1)(n-2)\binom{n-3}{k-3}\Bigr)\lambda-\Bigl(r(k-1)-(n-1)\binom{n-2}{k-2}\Bigr)^2=0.$
\end{theorem}

\begin{figure}[H]
\centering
\begin{tikzpicture}[scale=0.5]
   \draw (0,0) ellipse (1.2 and 3.5);
   \draw (3,0) ellipse (1.2 and 3.5);
   \node at (-.5,1.3) {$v_1$};
    \node at (0.4,1.3) {$v_2$};
    \node at (3,1.3) {$v_3$};
     \draw (1.5,1.5) ellipse (3.5 and 1.2);
     \draw (1.5,-1.5) ellipse (3.5 and 1.2);
     \node at (0,-1.7) {$v_4$};
    \node at (2.6,-1.7) {$v_5$};
    \node at (3.4,-1.7) {$v_6$};
\end{tikzpicture}

  \caption{$(3,2)$-regular hypergraph $\mathscr{G}^*$   }
  \label{fig:3}
\end{figure}
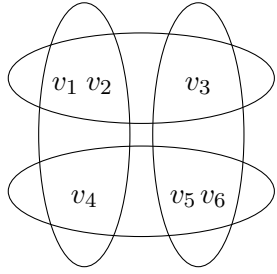

\begin{proof}
    The characteristic polynomial of $NNS(\mathscr{G}^*)$ is given by,
    \begin{align*}
      \det(\lambda I_n-NNS(\mathscr{G}^*))&=\det\begin{pmatrix}
   \lambda I_n- (n-2)\binom{n-3}{k-3}(J_n-I_n)& -\binom{n-2}{k-2}(J_n-I_n)+A(\mathscr{G}^*)\\
   \hspace{2.4cm}+(k-3)A(\mathscr{G}^*)& \\
	-\binom{n-2}{k-2}(J_n-I_n)+A(\mathscr{G}^*) &\lambda I_n
\end{pmatrix} \\
&=\det(\lambda I_n)\det\Bigl(\lambda I_n- (n-2)\binom{n-3}{k-3}(J_n-I_n)+(k-3)A(\mathscr{G}^*)\\&\hspace{6cm}-\displaystyle\frac{(A(\mathscr{G}^*)-\binom{n-2}{k-2}(J_n-I_n))^2}{\lambda}\Bigr)\\
&=\det\Bigl(\lambda^2 I_n- \lambda(n-2)\binom{n-3}{k-3}(J_n-I_n)+\lambda(k-3)A(\mathscr{G}^*)\\&\hspace{5.5cm}-(A(\mathscr{G}^*)-\binom{n-2}{k-2}(J_n-I_n))^2\Bigr).
    \end{align*}
Since  $A(\mathscr{G}^*)J_n=J_nA(\mathscr{G}^*)=r(k-1)J_n$. The characteristic polynomial of $NSS(\mathscr{G}^*)$ is, 
\begin{align*}
  P_{NNS(\mathscr{G}^*)}(\lambda)
  &=\lambda^2+\Bigl(r(k-1)(k-3)-(n-1)(n-2)\binom{n-3}{k-3}\Bigr)\lambda-\Bigl(r(k-1)-(n-1)\binom{n-2}{k-2}\Bigr)^2\\
  &\hspace{4cm}\prod_{i=2}^n \Bigl( \lambda^2+\Bigl((n-2)\binom{n-3}{k-3}+(k-3)\lambda_i\Bigr)\lambda-\Bigl(\lambda_i+\binom{n-2}{k-2}\Bigr)^2 \Bigr).
 \end{align*}
 Hence the result follows.
\end{proof}
\begin{corollary}\label{EnergyNNS}
     Let $\mathscr{G}^*$ be a $(k,r)-$regular hypergraph of order $n$ with eigenvalues $\lambda_1=r(k-1)\geq \lambda_2\geq \lambda_3\geq \ldots \geq \lambda_n$. If $k\geq 3$, then
     \begin{align*}
          \mathcal{E}(NNS(\mathscr{G}^*))&= \sqrt{\Bigl( r(k-1)(k-3)-(n-1)(n-2)\binom{n-3}{k-3}\Bigr)^2+4\Bigl(r(k-1)-(n-1)\binom{n-2}{k-2}\Bigr)^2}\\
          &\hspace{4cm}+\sum_{i=2}^n \sqrt{\bigl((k-3)\lambda_i+(n-2)\binom{n-3}{k-3}\bigr)^2+4\bigl(\lambda_i+\binom{n-2}{k-2}\bigr)^2}.
     \end{align*}
\end{corollary}
\begin{proof}
Take $X= r(k-1)(k-3)-(n-1)(n-2)\binom{n-3}{k-3},$ then 
    \begin{align*}
        \alpha_1&=\frac{-X+\sqrt{X^2+4\Bigl(r(k-1)-(n-1)\binom{n-2}{k-2}\Bigr)^2}}{2}\geq 0,\\
        \alpha_2&=\frac{-X-\sqrt{X^2+4\Bigl(r(k-1)-(n-1)\binom{n-2}{k-2}\Bigr)^2}}{2}\leq 0.
    \end{align*}
    Also,
    \begin{align*}
        \frac{-\bigl((k-3)\lambda_i+(n-2)\binom{n-3}{k-3}\bigr)+\sqrt{\bigl((k-3)\lambda_i+(n-2)\binom{n-3}{k-3}\bigr)^2+4\bigl(\lambda_i+\binom{n-2}{k-2}\bigr)}\bigr)^2}{2}   &\geq 0,\\
        \frac{-\bigl((k-3)\lambda_i+(n-2)\binom{n-3}{k-3}\bigr)-\sqrt{\bigl((k-3)\lambda_i+(n-2)\binom{n-3}{k-3}\bigr)^2+4\bigl(\lambda_i+\binom{n-2}{k-2}\bigr)}\bigr)^2}{2}&\leq 0.
    \end{align*}
    Hence, 
    \begin{equation*}
    \begin{split}
        |\alpha_1|+|\alpha_2|+\sum_{i=2}^n\bigg | \displaystyle \frac{-\bigl((k-3)\lambda_i+(n-2)\binom{n-3}{k-3}\bigr)\pm\sqrt{\bigl((k-3)\lambda_i+(n-2)\binom{n-3}{k-3}\bigr)^2+4\bigl(\lambda_i+\binom{n-2}{k-2}\bigr)}\bigr)^2}{2} \bigg | \\
        = \sqrt{\Bigl( r(k-1)(k-3)-(n-1)(n-2)\binom{n-3}{k-3}\Bigr)^2+4\Bigl(r(k-1)-(n-1)\binom{n-2}{k-2}\Bigr)^2}\\
          +\sum_{i=2}^n \sqrt{\bigl((k-3)\lambda_i+(n-2)\binom{n-3}{k-3}\bigr)^2+4\bigl(\lambda_i+\binom{n-2}{k-2}\bigr)^2}.
    \end{split}
    \end{equation*}
\end{proof}
Note that energy of complete $k$-uniform hypergraph on $k$-vertices, $\mathcal{E}(K_n^k)=2r(k-1)$. Applying Corollary \ref{EnergyNNS}, we get $\mathcal{E}(NNS(K^n_k))=2r(k-1).$ It is clear that for complete $k$-uniform hypergraph $K_n^k,$  $\mathcal{E}(K_n^k)=\mathcal{E}(NNS(K^k_n)).$\\
Following corollary shows that energy of $K_n^k$ act as a lower bound for the energy of both $NNS(K_n^k)$ and $NS(K_n^k).$
\begin{corollary}
    For a $k$-uniform hypergraph on $k$-vertices$(k\geq3)$, then
    $$\mathcal{E}(K_n^k)= \mathcal{E}(NNS(K_n^k)) < \mathcal{E}(NS(K_n^k)). $$
\end{corollary}
So it is natural to ask whether it is possible to compare the energy of $\mathscr{G}^*,\, NS(\mathscr{G}^*)$ and $NNS(\mathscr{G}^*)$. The following example illustrates that the inequality does not hold for all choices of $(k,r)$-regular hypergraph $\mathscr{G}^*$. Also, the example shows that the conjecture relating energy and nullity of graphs holds for $NS(G^*)$ and $NNS(G^*)$.
\begin{example}
   For  the hypergraph $\mathscr{G^*}$ in Figure \ref{fig:3}, $$\mathcal{E}(NNS(\mathscr{G}^*))=76.2998> \mathcal{E}(NS(\mathscr{G}^*))=24\sqrt{2}$$and
   $$\eta(NNS(\mathscr{G}^*))=0 <\eta(NS(\mathscr{G}^*))=2.$$
\end{example}
  We can have two immediate results on energy of $(k,r)$-regular hypergraph $\mathscr{G}^*, \,NS(\mathscr{G}^*)$ and $NNS(\mathscr{G}^*)$ from Corollaries \ref{EnergyNS} and \ref{EnergyNNS}. 
\begin{itemize}
    \item $\mathcal{E}(\mathscr{G}^*)<\mathcal{E}(NS_m(\mathscr{G}^*))$.
    \item $\mathcal{E}(\mathscr{G}^*)<\mathcal{E}(NNS(\mathscr{G}^*))$.
\end{itemize}
\section{Spectrum of Neighbourhood Splitting V-vertex and S-vertex Join}
\begin{definition}
    Let $\mathscr{G}_1^*$ and $\mathscr{G}_2^*$ be two $k$-uniform hypergraphs of order $n_1$ and $n_2$ respectively. Then the neighbourhood splitting $V$-vertex join of $\mathscr{G}_1^*$ and $\mathscr{G}_2^*$, denoted by $\mathscr{G}_1^*\veebar \mathscr{G}_2^*$ is the $k$-uniform hypergraph obtained from $NS(\mathscr{G}_1^*)$ and $\mathscr{G}_2^*$ by introducing new edges in such a way that each vertex of $V(\mathscr{G}_1^*)$ is a neighbour of every $(k-1)$ subset of the vertex set $V(\mathscr{G}_2^*)$.
\end{definition}
By arranging the vertices of $\mathscr{G}_1^*\veebar \mathscr{G}_2^*$ suitably, we obtain
\begin{equation*}
    A(\mathscr{G}_1^*\veebar \mathscr{G}_2^*)=\begin{bmatrix}
        (k-1)A(\mathscr{G}_1^*)&A(\mathscr{G}_1^*)&\displaystyle\binom{n_2-1}{k-2}J_{n_1,n_2}\\
        A(\mathscr{G}_1^*)&\textbf{0}_{n_1}&\textbf{0}_{n_1,n_2}\\
       \displaystyle \binom{n_2-1}{k-2}J_{n_2,n_1}&\textbf{0}_{n_2,n_1}&A(\mathscr{G}_2^*)+\displaystyle\binom{n_2-2}{k-3}n_1\left(J_{n_2}-I_{n_2}\right)
    \end{bmatrix}.
\end{equation*}
Here we obtain the characteristic polynomial of neighbourhood splitting V-vertex join of two $k-$uniform hypergraphs $\mathscr{G}_1^*$ and $\mathscr{G}_2^*$.
\begin{theorem}\label{thm4.2}
    Let $\mathscr{G}_i^*,i\in[1,2]$ be a $k$- uniform hypergraph on $n_i$ vertices and $\lambda_1^{(i)}\geq \lambda_2^{(i)} \geq \lambda_3^{(i)}\geq \ldots \geq \lambda_{n_i}^{(i)}$ are the adjacency eigenvalues of $A(\mathscr{G}_i^*)$. Then the characteristic polynomial of the neighbourhood splitting $V$-vertex join of the hypergraphs $\mathscr{G}_1^*$ and $\mathscr{G}_2^*$ is
     \begin{align*}
      P_{A(\mathscr{G}_1^* \veebar \mathscr{G}_2^*)}(\lambda)=&\left(\lambda+n_1(1-n_2)\binom{n_2-2}{k-3}-\lambda_1^{(2)}\right)\prod_{i=2}^{n_2}\left(\lambda+\binom{n_2-2}{k-3}n_1-\lambda_i^{(2)} \right)\\ &\prod_{i=1}^{n_1} \left(\lambda^2-(k-1)\lambda_i^{\scriptscriptstyle(1)}\lambda-{\lambda_i^{\scriptscriptstyle(1)}}^2\right)\left( 1-\displaystyle\binom{n_2-1}{k-2}^2 \rchi_{R}(\lambda)\rchi_{R'}(\lambda) \right),
 \end{align*} 
 where $\displaystyle R=A(\mathscr{G}_2^*)+\displaystyle\binom{n_2-2}{k-3}n_1\left(J_{n_2}-I_{n_2}\right)$ and $R'=(k-1)A(\mathscr{G}_1^*)+\frac{A(\mathscr{G}_1^*)^2}{\lambda}$.
\end{theorem}
\begin{proof}
    The characteristic polynomial of $A(\mathscr{G}_1^* \veebar \mathscr{G}_2^*)$ is
    \begin{equation*}
        P_{A(\mathscr{G}_1^* \veebar \mathscr{G}_2^*)}(\lambda)=\begin{pmatrix}
           \lambda I_{n_1}-(k-1)A(\mathscr{G}_1^*)&-A(\mathscr{G}_1^*)&-\displaystyle\binom{n_2-1}{k-2}J_{n_1,n_2}\\
        -A(\mathscr{G}_1^*)&\lambda I_{n_1}&\textbf{0}_{n_1,n_2}\\
       -\displaystyle \binom{n_2-1}{k-2}J_{n_2,n_1}&\textbf{0}_{n_2,n_1}&\lambda I_{n_2}-A(\mathscr{G}_2^*)-\displaystyle\binom{n_2-2}{k-3}n_1\left(J_{n_2}-I_{n_2}\right)
        \end{pmatrix}.
    \end{equation*}
    Using Lemma \ref{detblock}, we get
   \begin{equation}\label{eqn2}
        P_{A(\mathscr{G}_1^* \veebar \mathscr{G}_2^*)}(\lambda)=\det (\lambda I_{n_2}-R) \det (S),
    \end{equation}
    where
    $$S=\begin{bmatrix}
         \lambda I_{n_1}-(k-1)A(\mathscr{G}_1^*)&-A(\mathscr{G}_1^*)\\
        -A(\mathscr{G}_1^*)&\lambda I_{n_1}
    \end{bmatrix}-\binom{n_2-1}{k-2}^2\begin{bmatrix}
        J_{n_1,n_2}\\
        \textbf{0}_{n_1,n_2}
    \end{bmatrix}(\lambda I_{n_2}-R)^{-1}\begin{bmatrix}
         J_{n_2,n_1}&\textbf{0}_{n_2,n_1}
    \end{bmatrix}$$
    and $\displaystyle R=A(\mathscr{G}_2^*)+\displaystyle\binom{n_2-2}{k-3}n_1\left(J_{n_2}-I_{n_2}\right).$ Applying Definition \ref{coronaldefn}, we get
    \begin{align*}
       S&=\begin{bmatrix}
         \lambda I_{n_1}-(k-1)A(\mathscr{G}_1^*)&-A(\mathscr{G}_1^*)\\
        -A(\mathscr{G}_1^*)&\lambda I_{n_1}
    \end{bmatrix}- \displaystyle\binom{n_2-1}{k-2}^2 \begin{bmatrix}
        I_{n_1,n_2}\\
        \textbf{0}_{n_1,n_2}
    \end{bmatrix} J_{n_2,1} \rchi_{R}(\lambda)J_{1,n_2} \begin{bmatrix}
         I_{n_2,n_1}&\textbf{0}_{n_2,n_1}
    \end{bmatrix} \\
    &=\begin{bmatrix}
         \lambda I_{n_1}-(k-1)A(\mathscr{G}_1^*)&-A(\mathscr{G}_1^*)\\
        -A(\mathscr{G}_1^*)&\lambda I_{n_1}
    \end{bmatrix}- \displaystyle\binom{n_2-1}{k-2}^2 \rchi_{R} \begin{bmatrix}
        J_{n_1} & \textbf{0}_{n_1}\\
        \textbf{0}_{n_1} &\textbf{0}_{n_1}
    \end{bmatrix} \\
    &=\begin{bmatrix}
         \lambda I_{n_1}-(k-1)A(\mathscr{G}_1^*)-\displaystyle\binom{n_2-1}{k-2}^2 \rchi_{R} J_{n_1}&-A(\mathscr{G}_1^*)\\
        -A(\mathscr{G}_1^*)&\lambda I_{n_1}
    \end{bmatrix}.
    \end{align*}
    Again by Lemma \ref{detblock}, we obtain
    \begin{align*}
        \det(S)=\det(\lambda I_{n_1})\det\left( \lambda I_{n_1}-(k-1)A(\mathscr{G}_1^*)-\displaystyle\binom{n_2-1}{k-2}^2 \rchi_{R} J_{n_1}-\frac{A(\mathscr{G}_1^*)^2}{\lambda}\right).
    \end{align*}
   Take $R'=(k-1)A(\mathscr{G}_1^*)+\frac{A(\mathscr{G}_1^*)^2}{\lambda}$. From Lemma \ref{detMJ}, we have
   \begin{align*}
       \det(S)&=\det(\lambda I_{n_1})\left( \det(\lambda I_{n_1}-R')-\displaystyle\binom{n_2-1}{k-2}^2 \rchi_{R}(\lambda) J_{1,n_1} adj(\lambda I_{n_1}-R')J_{n_1,1}\right)\\
              &=\det(\lambda I_{n_1})\det(\lambda I_{n_1}-R') \left( 1-\displaystyle\binom{n_2-1}{k-2}^2 \rchi_{R}(\lambda) J_{1,n_1} (\lambda I_{n_1}-R')^{-1}J_{n_1,1}\right)\\
              &=\lambda^{n_1}\det(\lambda I_{n_1}-R') \left( 1-\displaystyle\binom{n_2-1}{k-2}^2 \rchi_{R}(\lambda)\rchi_{R'}(\lambda) \right).
   \end{align*}
 Since $\det(\lambda I_{n_1}-R')=\displaystyle\frac{1}{\lambda^{n_1}} \displaystyle\prod_{i=1}^{n_1} \lambda^2-(k-1)\lambda_i^{\scriptscriptstyle(1)}\lambda-{\lambda_i^{\scriptscriptstyle(1)}}^2$, we have
 \begin{equation}\label{eqn3}
     \det(S)=\prod_{i=1}^{n_1} \left(\lambda^2-(k-1)\lambda_i^{\scriptscriptstyle(1)}\lambda-{\lambda_i^{\scriptscriptstyle(1)}}^2\right)\left( 1-\displaystyle\binom{n_2-1}{k-2}^2 \rchi_{R}(\lambda)\rchi_{R'}(\lambda) \right).
 \end{equation}
Now consider
 \begin{align}\label{eqn4}
     \det (\lambda I_{n_2}-R) 
     & =\det\left((\lambda+\displaystyle\binom{n_2-2}{k-3}n_1)I_{n_2}-A(\mathscr{G}_2^*)-\displaystyle\binom{n_2-2}{k-3}n_1 J_{n_2}\right)\notag \\
    &=\left(\lambda+n_1(1-n_2)\binom{n_2-2}{k-3}-\lambda_1^{(2)}\right)\prod_{i=2}^{n_2}\left(\lambda+\binom{n_2-2}{k-3}n_1-\lambda_i^{(2)} \right).
 \end{align}
 From (\ref{eqn2}),(\ref{eqn3}) and (\ref{eqn4}), we get
 \begin{align*}
      P_{A(\mathscr{G}_1^* \veebar \mathscr{G}_2^*)}(\lambda)=&\left(\lambda+n_1(1-n_2)\binom{n_2-2}{k-3}-\lambda_1^{(2)}\right)\prod_{i=2}^{n_2}\left(\lambda+\binom{n_2-2}{k-3}n_1-\lambda_i^{(2)} \right)\\ &\hspace{.5cm}\prod_{i=1}^{n_1} \left(\lambda^2-(k-1)\lambda_i^{\scriptscriptstyle(1)}\lambda-{\lambda_i^{\scriptscriptstyle(1)}}^2\right)\left( 1-\displaystyle\binom{n_2-1}{k-2}^2 \rchi_{R}(\lambda)\rchi_{R'}(\lambda) \right).
 \end{align*} 
 Hence the theorem.
\end{proof}
\begin{theorem}
      Let $\mathscr{G}_i^*,i\in[1,2]$ be a $(k,r_i)$- regular hypergraph on $n_i$ vertices and $\lambda_1^{(i)}=r_i(k-1)\geq \lambda_2^{(i)} \geq \lambda_3^{(i)}\geq \ldots \geq \lambda_{n_i}^{(i)}$ are the adjacency eigenvalues of $A(\mathscr{G}_i^*)$. Then the adjacency spectrum of the neighbourhood splitting $V$-vertex join of the hypergraphs $\mathscr{G}_1^*$ and $\mathscr{G}_2^*$ is,
       $$\sigma_A({\mathscr{G}_1^*\veebar \mathscr{G}_2^*})=
    \begin{pmatrix}
    \lambda_i^{(2)}-n_1\binom{n_2-2}{k-3} &\lambda_i^{(1)}\frac{(k-1+\sqrt{(k-1)^2+4})}{2} &\lambda_i^{(1)}\frac{(k-1-\sqrt{(k-1)^2+4)}}{2} &\alpha_1&\alpha_2&\alpha_3\\
             1&1&1&1&1&1
    \end{pmatrix},$$
    where $\alpha_1,\alpha_2$ and $\alpha_3$ are the roots of the equation $\lambda^3-(a+r_1(k-1)^2)\lambda^2+\Bigl(r_1(k-1)^2(a-r_1)-n_1n_2\binom{n_2-1}{k-2}^2\Bigr)\lambda+a r_1^2(k-1)^2=0, \,a=\displaystyle r_2(k-1)+\binom{n_2-2}{k-3}n_1(n_2-1).$
\end{theorem}
\begin{proof}
    From Theorem \ref{thm4.2}, we get
 \begin{align*}
      P_{A(\mathscr{G}_1^* \veebar \mathscr{G}_2^*)}(\lambda)=&\left(\lambda+n_1(1-n_2)\binom{n_2-2}{k-3}-r(k-1)\right)\left(\lambda^2-r_1(k-1)^2\lambda-r_1^2(k-1)^2\right)\\ & \prod_{i=2}^{n_2}\left(\lambda+\binom{n_2-2}{k-3}n_1-\lambda_i^{(2)} \right) \prod_{i=2}^{n_1} \left(\lambda^2-(k-1)\lambda_i^{\scriptscriptstyle(1)}\lambda-{\lambda_i^{\scriptscriptstyle(1)}}^2\right)\\ &\left( 1-\displaystyle\binom{n_2-1}{k-2}^2 \rchi_{R}(\lambda)\rchi_{R'}(\lambda) \right),
 \end{align*} 
 where $\displaystyle  \rchi_{R}(\lambda)=\frac{n_2}{\lambda-r_2(k-1)-\binom{n_2-2}{k-3}n_1(n_2-1)}$ and $\rchi_{R'}(\lambda)=\displaystyle\frac{n_1\lambda}{\lambda^2-r_1(k-1)^2\lambda+r_1^2(k-1)^2}$. On simplification, we get the desired result.
\end{proof}
\begin{corollary}
    Let $\mathscr{G}_1^*$ and $\mathcal{H}_1^*$ be $(k,r_1)$- regular hypergraphs on $n_1$ vertices and $\mathscr{G}_2^*$ be a $(k,r_2)$- regular hypergraph on $n_2$ vertices. If $\mathscr{G}_1^*$ and $\mathcal{H}_1^*$ are non-isomorphic cospectral, then $\mathscr{G}_1^*\veebar \mathscr{G}_2^*$ and $\mathcal{H}_1^*\veebar \mathscr{G}_2^*$, and $\mathscr{G}_2^*\veebar \mathscr{G}_1^*$ and $\mathscr{G}_2^*\veebar \mathcal{H}_1^*$ are non-isomorphic cospectral. 
\end{corollary}
\begin{corollary}
    Let $\mathscr{G}_i^*,i\in[1,2]$ be a $(k,r_i)$- regular hypergraph on $n_i$ vertices. If $\mathscr{G}_1^*$ be a singular hypergraph, then $\mathscr{G}_1^*\veebar \mathscr{G}_2^*$ is also a singular hypergraph.
\end{corollary}
\begin{definition}
    Let $\mathscr{G}_1^*$ and $\mathscr{G}_2^*$ be two $k$-uniform hypergraphs of order $n_1$ and $n_2$ respectively. Then the neighbourhood splitting $S$-vertex join of $\mathscr{G}_1^*$ and $\mathscr{G}_2^*$, denoted by $\mathscr{G}_1^*\barwedge \mathscr{G}_2^*$ is the $k$-uniform hypergraph obtained from $NS(\mathscr{G}_1^*)$ and $\mathscr{G}_2^*$ by introducing new edges in such a way that each vertex of $S(\mathscr{G}_1^*)$ is a neighbour of every $(k-1)$ subset of the vertex set $V(\mathscr{G}_2^*)$.
\end{definition}
\noindent By arranging the vertices of $\mathscr{G}_1^*\barwedge \mathscr{G}_2^*$ suitably, we obtain
\begin{equation*}
    A(\mathscr{G}_1^*\barwedge \mathscr{G}_2^*)=\begin{bmatrix}
        (k-1)A(\mathscr{G}_1^*)&A(\mathscr{G}_1^*)&\textbf{0}_{n_1,n_2}\\
        A(\mathscr{G}_1^*)&\textbf{0}_{n_1}&\displaystyle\binom{n_2-1}{k-2}J_{n_1,n_2}\\
       \textbf{0}_{n_2,n_1}&\displaystyle \binom{n_2-1}{k-2}J_{n_2,n_1}&A(\mathscr{G}_2^*)+\displaystyle\binom{n_2-2}{k-3}n_1\left(J_{n_2}-I_{n_2}\right)
    \end{bmatrix}.
\end{equation*}
\begin{theorem}
      Let $\mathscr{G}_i^*$ be a $(k,r_i)$- regular hypergraph on $n_i$ vertices and $\lambda_1^{(i)}=r_i(k-1)\geq \lambda_2^{(i)} \geq \lambda_3^{(i)}\geq \ldots \geq \lambda_{n_i}^{(i)}$ are the adjacency eigenvalues of $A(\mathscr{G}_i^*)$. Then the characteristic polynomial of the neighbourhood splitting $S$-vertex join of the hypergraphs $\mathscr{G}_1^*$ and $\mathscr{G}_2^*$ is,
    \begin{equation*}
    \begin{split}
          P_{A(\mathscr{G}_1^* \barwedge \mathscr{G}_2^*)}(\lambda) =
          = &\biggl(\Bigl(\lambda-r_2(k-1)-n_1(n_2-1)\binom{n_2-2}{k-3}\Bigr)\Bigl(\lambda^2-r_1(k-1)^2\lambda-r_1^2(k-1)^2\Bigr)\\
         &\hspace{1pt}-n_1n_2\binom{n_2-1}{k-2}^2\Bigl(\lambda-r_1(k-1)^2\Bigr)\biggr)\,\prod_{j=2}^{n_2}\left(\lambda+n_1\binom{n_2-2}{k-3}-\lambda_i^{{\scriptscriptstyle(2)}}\right)\\ &\hspace{1pt}\prod_{j=2}^{n_1}\left(\lambda^2-(k-1)\lambda_i^{{\scriptscriptstyle(1)}}\lambda-\lambda_i^{{\scriptscriptstyle(1)}^2}\right).
    \end{split}
    \end{equation*}
    \end{theorem}
    \begin{proof}
          The characteristic polynomial of a matrix $A(\mathscr{G}_1^* \barwedge \mathscr{G}_2^*)$ is
    \begin{equation*}
        P_{A(\mathscr{G}_1^* \barwedge \mathscr{G}_2^*)}(\lambda)=\begin{pmatrix}
           \lambda I_{n_1}-(k-1)A(\mathscr{G}_1^*)&-A(\mathscr{G}_1^*)&\textbf{0}_{n_1,n_2}\\
        -A(\mathscr{G}_1^*)&\lambda I_{n_1}&-\displaystyle\binom{n_2-1}{k-2}J_{n_1,n_2}\\
      \textbf{0}_{n_2,n_1}& -\displaystyle \binom{n_2-1}{k-2}J_{n_2,n_1}&\substack{\lambda I_{n_2}-A(\mathscr{G}_2^*)\\-\binom{n_2-2}{k-3}n_1\left(J_{n_2}-I_{n_2}\right)}\\

        \end{pmatrix}.
    \end{equation*}
    From Lemma \ref{detblock}, we get
   \begin{equation}\label{eqn5}
        P_{A(\mathscr{G}_1^* \veebar \mathscr{G}_2^*)}(\lambda)=\det (\lambda I_{n_2}-R) \det (S),
    \end{equation}
    where
    $$S=\begin{bmatrix}
         \lambda I_{n_1}-(k-1)A(\mathscr{G}_1^*)&-A(\mathscr{G}_1^*)\\
        -A(\mathscr{G}_1^*)&\lambda I_{n_1}
    \end{bmatrix}-\binom{n_2-1}{k-2}^2\begin{bmatrix}
        \textbf{0}_{n_1,n_2}\\
         J_{n_1,n_2}
    \end{bmatrix}(\lambda I_{n_2}-R)^{-1}\begin{bmatrix}
         \textbf{0}_{n_2,n_1}&J_{n_2,n_1}
    \end{bmatrix}$$
    and $\displaystyle R=A(\mathscr{G}_2^*)+\displaystyle\binom{n_2-2}{k-3}n_1\left(J_{n_2}-I_{n_2}\right).$ Applying Definition \ref{coronaldefn}, we get
    \begin{align*}
       S&=\begin{bmatrix}
         \lambda I_{n_1}-(k-1)A(\mathscr{G}_1^*)&-A(\mathscr{G}_1^*)\\
        -A(\mathscr{G}_1^*)&\lambda I_{n_1}
    \end{bmatrix}- \displaystyle\binom{n_2-1}{k-2}^2 
    \begin{bmatrix}
        \textbf{0}_{n_1,n_2}\\
           I_{n_1,n_2}
    \end{bmatrix} J_{n_2,1} \rchi_{R}(\lambda)J_{1,n_2} \begin{bmatrix}
         \textbf{0}_{n_2,n_1}&I_{n_2,n_1}
    \end{bmatrix} \\
    &=\begin{bmatrix}
         \lambda I_{n_1}-(k-1)A(\mathscr{G}_1^*)&-A(\mathscr{G}_1^*)\\
        -A(\mathscr{G}_1^*)&\lambda I_{n_1}
    \end{bmatrix}- \displaystyle\binom{n_2-1}{k-2}^2 \rchi_{R}(\lambda) \begin{bmatrix}
        \textbf{0}_{n_1} & \textbf{0}_{n_1}\\
        \textbf{0}_{n_1} &J_{n_1}
    \end{bmatrix} \\
    &=\begin{bmatrix}
         \lambda I_{n_1}-(k-1)A(\mathscr{G}_1^*)&-A(\mathscr{G}_1^*)\\
        -A(\mathscr{G}_1^*)&\lambda I_{n_1}-\displaystyle\binom{n_2-1}{k-2}^2 \rchi_{R}(\lambda) J_{n_1}
    \end{bmatrix}.
    \end{align*} 
      Therefore, we have
    \begin{align*}
    \begin{split}
         \det(S)=\det\Bigl(\lambda I_{n_1}&-\displaystyle\binom{n_2-1}{k-2}^2 \rchi_{R}(\lambda) J_{n_1}\Bigr)\det\biggl( \lambda I_{n_1}-(k-1)A(\mathscr{G}_1^*)\\
         &\hspace{2pt}-A(\mathscr{G}_1^*)\Bigl( \lambda I_{n_1}-\displaystyle\binom{n_2-1}{k-2}^2 \rchi_{R}(\lambda) J_{n_1}\Bigr)^{-1}A(\mathscr{G}_1^*)\biggr).
    \end{split}
    \end{align*}
    By Lemma \ref{invIJ} and on simplification, we get
    \begin{equation*}
        \begin{split}
         \det(S)=\lambda^{n_1-1}\Bigl(\lambda-n_1\displaystyle\binom{n_2-1}{k-2}^2 \rchi_{R}(\lambda)\Bigr)\det\biggl( \lambda I_{n_1}-R'-\displaystyle\frac{r_1^2(k-1)^2\binom{n_2-1}{k-2}^2  \rchi_{R}(\lambda)}{\lambda (\lambda-n_1 \binom{n_2-1}{k-2}^2  \rchi_{R}(\lambda))}J_{n_1}\biggr),
    \end{split}
    \end{equation*}
  where $R'=(k-1)A(\mathscr{G}_1^*)+\displaystyle\frac{1}{\lambda}A^2(\mathscr{G}_1^*)$. From Lemma \ref{detMJ}, we have
   \begin{align*}
       \det(S)&=\lambda^{n_1-1}\Bigl(\lambda-n_1\displaystyle\binom{n_2-1}{k-2}^2 \rchi_{R}(\lambda)\Bigr)\det(\lambda I_{n_1}-R')\biggl( 1-\displaystyle\frac{r_1^2(k-1)^2\binom{n_2-1}{k-2}^2  \rchi_{R}(\lambda)\rchi_{R'}(\lambda)}{\lambda (\lambda-n_1 \binom{n_2-1}{k-2}^2  \rchi_{R}(\lambda))}\biggr)\\
       &=\lambda^{n_1-2}\det(\lambda I_{n_1}-R')\biggl(\lambda\Bigl(\lambda-n_1\displaystyle\binom{n_2-1}{k-2}^2 \rchi_{R}(\lambda)\Bigr)\\&\hspace{6cm}-r_1^2(k-1)^2\binom{n_2-1}{k-2}^2  \rchi_{R}(\lambda)\rchi_{R'}(\lambda)\biggr).
   \end{align*}
 Since $\det(\lambda I_{n_1}-R')=\displaystyle\frac{1}{\lambda^{n_1}} \displaystyle\prod_{i=1}^{n_1} \lambda^2-(k-1)\lambda_i^{\scriptscriptstyle(1)}\lambda-{\lambda_i^{\scriptscriptstyle(1)}}^2$, we have
 \begin{equation}\label{eqn6}
 \begin{split}
       \det(S)=\frac{1}{\lambda^2}\Bigl(\lambda^2-n_1\displaystyle\binom{n_2-1}{k-2}^2 \rchi_{R}(\lambda)\lambda-r_1^2&(k-1)^2\binom{n_2-1}{k-2}^2  \rchi_{R}(\lambda)\rchi_{R'}(\lambda)\Bigr)\\ &\hspace{10pt}\prod_{i=1}^{n_1} \left(\lambda^2-(k-1)\lambda_i^{\scriptscriptstyle(1)}\lambda-{\lambda_i^{\scriptscriptstyle(1)}}^2\right),
\end{split}
\end{equation}
where $\rchi_{R}(\lambda)=\displaystyle\frac{n_2}{\lambda-r_2(k-1)-\binom{n_2-2}{k-3}n_1(n_2-1)} \text{ and }\rchi_{R'}(\lambda)=\displaystyle\frac{n_1}{\lambda-r_1(k-1)^2-\frac{1}{\lambda}r_1^2(k-1)^2}.$\\
Now consider
\begin{align}\label{eqn7}
    \det (\lambda I_{n_2}-R)&=\det \Bigl(\lambda I_{n_2}-A(\mathscr{G}_2^*)-\binom{n_2-2}{k-3}n_1(J_{n_2}-I_{n_2})\Bigr)\notag\\
    &=\Bigl(\lambda-r_2(k-1)-\binom{n_2-2}{k-3}n_1(n_2-1)\Bigr)\prod_{i=2}^{n_2}\Bigl(\lambda+n_1\binom{n_2-2}{k-3}-\lambda_i^{\scriptscriptstyle(2)}\Bigr).
\end{align}
From (\ref{eqn5}),(\ref{eqn6}) and (\ref{eqn7}) we get the desired result.
    \end{proof}
 \begin{corollary}
     Let $\mathscr{G}_i^*$ be a $(k,r_i)$- regular hypergraph on $n_i$ vertices and $\lambda_1^{(i)}=r_i(k-1)\geq \lambda_2^{(i)} \geq \lambda_3^{(i)}\geq \ldots \geq \lambda_{n_i}^{(i)}$ are the adjacency eigenvalues of $A(\mathscr{G}_i^*)$. Then the spectrum of $\mathscr{G}_1^* \barwedge \mathscr{G}_2^*$ is,
    \begin{equation*}
           \sigma_{A}(\mathscr{G}_1^* \barwedge \mathscr{G}_2^*) 
           =\begin{pmatrix}
              -n_1\binom{n_2-2}{k-3}+\lambda_i^{{\scriptscriptstyle(2)}}&\frac{\lambda_i^{{\scriptscriptstyle(1)}^2}((k-1)+\sqrt{(k-1)^2+4})}{2}&\frac{\lambda_i^{{\scriptscriptstyle(1)}^2}((k-1)-\sqrt{(k-1)^2+4})}{2}&\alpha_1&\alpha_2&\alpha_3\\
              1& 1&1&1&1&1
             \end{pmatrix},
    \end{equation*}
    where $i\in[2,n_i]$, and $\alpha_1,\alpha_2$ and $\alpha_3$ are the roots of the equation $\Bigl(\lambda-r_2(k-1)-n_1(n_2-1)\binom{n_2-2}{k-3}\Bigr)\Bigl(\lambda^2-r_1(k-1)^2\lambda-r_1^2(k-1)^2\Bigr)-n_1n_2\binom{n_2-1}{k-2}^2\Bigl(\lambda-r_1(k-1)^2\Bigr).$
 \end{corollary}
 
\begin{corollary}
    Let $\mathscr{G}_1^*$ and $\mathcal{H}_1^*$ be $(k,r_1)$- regular hypergraphs of order $n_1$ and $\mathscr{G}_2^*$ be a $(k,r_2)$- regular hypergraph of order $n_2$. If $\mathscr{G}_1^*$ and $\mathcal{H}_1^*$ are non-isomorphic cospectral hypergraphs, then $\mathscr{G}_1^* \barwedge \mathscr{G}_2^*$ and $\mathcal{H}_1^*\barwedge \mathscr{G}_2^*$, as well as $\mathscr{G}_2^*\barwedge \mathscr{G}_1^*$ and $\mathscr{G}_2^*\barwedge \mathcal{H}_1^*$ are also non-isomorphic cospectral hypergraphs. 
\end{corollary}
\begin{corollary}
    Let $\mathscr{G}_i^*,i\in[1,2]$ be a $(k,r_i)$- regular hypergraph on $n_i$ vertices. If $\mathscr{G}_1^*$ be a singular hypergraph, then $\mathscr{G}_1^*\barwedge \mathscr{G}_2^*$ is also a singular hypergraph.
\end{corollary}
\section{Conclusion}   
 In this paper, the spectra of neighbourhood m-splitting and non-neighbourhood splitting hypergraphs, as well as their respective energies are determined. In addition, the spectrum of neighbourhood splitting V-vertex and S-vertex join of hypergraphs are estimated. As an application, an infinite families of singular graphs and infinite pairs of non-isomorphic cospectral hypergraphs are constructed. 

\section{Declarations}
 The authors declare that there is no conflict of interest.\\
\bibliographystyle{plain}
\bibliography{ref}
\end{document}